\documentclass[a4paper,12pt]{article}
\usepackage{amsthm}
\usepackage{amsfonts}
\usepackage{amsmath}
\usepackage{amssymb}
\usepackage{diagrams}
\usepackage{hyperref}
\usepackage[pdftex]{graphicx}
\usepackage{chngcntr}
\usepackage{thmtools}
\usepackage{thm-restate}
\counterwithin{table}{section}

\newtheorem{theorem}{Theorem}[section]

\newtheorem{proposition}[theorem]{Proposition}

\theoremstyle{definition}

\newtheorem*{definition*}{Definition}

\theoremstyle{remark}

\numberwithin{equation}{section}

\newcommand{\Z}{\mathbb{Z}}

\newcommand{\Q}{\mathbb{Q}}

\newcommand{\diag}{\operatorname{diag}}
\newcommand{\bq}{/\!\!/}

\begin{document}

\title{Rationally 4-periodic biquotients}

\author{Jason DeVito}

\date{\vspace{-5ex}}

\maketitle

\abstract{An $n$-dimensional manifold $M$ is said to be rationally $4$-periodic if there is an element $e\in H^4(M;\Q)$ with the property that cupping with $e$, $\cdot \cup e:H^\ast(M;\Q)\rightarrow H^{\ast + 4}(M;\Q)$ is injective for $0< \ast \leq \dim M-4$ and surjective when $0\leq \ast < \dim M-4$.  We classify all compact simply connected biquotients which are rationally $4$-periodic.  In addition, we show that if a simply connected rationally elliptic CW-complex $X$ of dimension at least $6$ is rationally $4$-periodic, then the cohomology ring is either singly generated, or $X$ is rationally homotopy equivalent to $S^2\times \mathbb{H}P^n$, $S^3\times \mathbb{H}P^n$, or $S^3\times S^3$.}

\section{Introduction}

A biquotient is any manifold which is diffeomorphic to the quotient of a Riemannian homogeneous space by a free isometric action.  Biquotients are important in the study of positively curved manifolds.  In fact, with the exception of the positively curved manifold found independently by Dearricott \cite{De} and by Grove, Verdiani, and Ziller \cite{GVZ}, all known examples of manifolds admitting positive curvature are diffeomorphic to biquotients \cite{Ber,AW,Wa,Es1,Es2,Baz1,GVZ,De}.  Further, all known examples of manifolds admitting quasi or almost positive curvature are diffeomorphic to biquotients \cite{DDRW,Di1,EK,GrMe1,Ke1,Ke2,KT,PW2,Ta1,W,Wi}.

Recall that an $n$-dimensional manifold $M$ is said to be $k$-periodic with coefficient ring $R$ if there is an element $e\in H^k(M;R)$ with the property that cupping with $e, \cup e:H^\ast(M)\rightarrow H^{\ast + k}(M)$ is injective when $0<\ast\leq n-k$ and surjective when $0 \leq \ast < n-k$.  Due to Wilking's connectedness lemma \cite{Wi2}, periodicity appears when $M$ admits a metric of positive sectional curvature with large symmetry group.  For example, Wilking \cite{Wi2} has shown that if $M^n$ with $n\geq 6000$ is positively curved and admits an effective isometric $T^d$ action for $d\geq \frac{1}{6}n+1$, then $M$ is $4$-periodic with respect to any field coefficients.  Further, if the period k is small compared to the dimension of the manifold, Kennard \cite{Ken1} shows that k-periodic integral cohomology implies $4$-periodic rational cohomology. As an application, in \cite{Ken2}, Kennard shows that a simply connected closed manifold with positive sectional curvature and large symmetry rank has $4$-periodic rational cohomology in small degrees.

If a manifold $M^n$ has a $k$-periodic rational cohomology ring with $k=1$ or $2$, or if $k=4$ and the second Betti number vanishes, then it is easy to see $H^\ast(M;\Q)$ is isomorphic to that of a compact rank one symmetric space other than $\mathbb{O}P^2$.  So, the cohomology ring of such examples are generated by a single element.  Simply connected biquotients with singly generated rational cohomology rings have been classified by Kapovitch and Ziller \cite{KZ}.  In particular, every simply connected biquotient with singly generated cohomology is $4$-periodic, with the exception of $\mathbb{O}P^2$.  Thus, for biquotients, rational $4$-periodicity may be viewed as a generalization of having singly generated cohomology.

We note that if $M$ has dimension at most $4$, then it is vacuously $4$-periodic,  Similarly if $M$ is $5$ dimensional and simply connected, then it is $4$-periodic.  Compact simply connected biquotients of dimension at most $5$ were classified in \cite{To1,Pa1,DeV1}  Thus, we will assume the dimension of our biquotients to be at least $6$.

Our first theorem applies to a larger class of spaces, rationally elliptic spaces.  Recall that a simply connected CW complex $X$ is called rationally elliptic if the rational vector spaces $H^\ast(X;\Q)$ and $\bigoplus \pi_{\ast}(X)\otimes\Q$ are finite dimensional.  We use the notation $X\simeq_\Q Y$ to indicate that $X$ and $Y$ have the same rational homotopy type.

\begin{restatable}{theorem}{mainone}
\label{rattype}Suppose $X$ is an $n$-dimensional simply connected $CW$-complex with $n\geq 6$ which is rationally elliptic.  Further, assume $H^\ast(X;\mathbb{Q})$ is $4$-periodic but not singly generated.  Then precisely one of the following occurs.

(a)  $X\simeq_\Q\mathbb{H}P^m\times S^3$

(b)  $X\simeq_\Q\mathbb{H}P^m\times S^2$

(c)  $X\simeq_\Q S^3\times S^3$.

\end{restatable}

Under the assumption that $X$ is diffeomorphic to a biquotient and $X\simeq_\Q S^3\times S^3$, the author showed in \cite[Corollary 3.5]{DeV2} that $X$ is diffeomorphic to $S^3\times S^3$.  Using Kapovitch and Ziller's \cite{KZ} classification of biquotients with singly generated cohomology, we obtain a similar characterization in case (a), when $X\simeq_\Q \mathbb{H}P^m\times S^3$ is a biquotient.

\begin{theorem}\label{bclass}  Assume $M$ is a simply connected manifold diffeomorphic to a biquotient.  If $M\simeq_\Q \mathbb{H}P^m\times S^3$, then $M$ is diffeomorphic to exactly one of the following:

(1)  the total space of one of two (respectively three) $S^3$ bundles over $\mathbb{H}P^m$ for $m\geq 2$ (respectively $m=1$),

(2)  if $m=2$, the total space of one of three $S^3$ bundles over $\mathbf{G}_2/SO(4)$ , or

(3)  if $m\geq 3$ is odd, the total space of one of two $S^3$ bundles over the biquotient $\Delta SU(2)\backslash SO(2m+3)/SO(2m+1)$.

\end{theorem}

In particular, in dimension $7$ there are $3$ examples, in dimension $11$ there are $5$, in higher dimensions of the form $8m+3$ there are $2$ examples, and in higher dimensions of the form $8m-1$ there are $4$, all distinct up to diffeomorphism.  In (3), if we allow $m = 1$, then these biquotients are diffeomorphic to examples in $(1)$.

Unfortunately, we can only partially classify biquotients with the rational homotopy type of $\mathbb{H}P^m\times S^2$.

\begin{theorem}\label{cclass} Assume $M$ is a simply connected manifold diffeomorphic to a biquotient.  If $M\simeq_\Q \mathbb{H}P^m\times S^2$,  then $M$ is diffeomorphic to the exactly one of the following:

(1)   The total space of one of the two linear bundles over $S^2$ with fiber $\mathbb{H}P^m$ for $m\geq 1$.

(2)  if $m=2$, $[\mathbf{G}_2/SO(4)] \times S^2$

(3)  if $m\geq 3$ is odd, the total space of some bundle over $S^2$ with fiber the biquotient $\Delta SU(2)\backslash SO(2m+3)/SO(2m+1)$.

Further, in case (3), only finitely many diffeomorphism types arise for each $m$.

\end{theorem}

This paper is organized as follows.  In Section \ref{reduct}, we recall the basics of rational homotopy theory and biquotients, proving Theorem \ref{rattype}.  In Section \ref{S3}, we begin with a structure theorem for biquotients which are rationally $\mathbb{H}P^m\times S^3$, see Theorem \ref{hpms3}.  We then classify biquotients having this structure, proving Theorem \ref{bclass}.  In Section \ref{S2}, we prove an analogous structure theorem, Theorem \ref{hpms3}, and use it to prove Theorem \ref{cclass}.

\

We would like to thank Manuel Amann and Lee Kennard for suggesting this problem, as well as for several stimulating discussions.

\section{Background}\label{reduct}

In this section, we recall the necessary background regarding rational homotopy theory and biquotients.

\subsection{Rational homotopy theory}\label{rht}

In this section, we prove Theorem \ref{rattype} after introducing the relevant notions from rational homotopy theory \cite{FeOpTa,FeHaTh}.  We use the notation $\pi_\ast(X)_\Q$ as shorthand for the rational homotopy group $\pi_\ast(X)\otimes \mathbb{Q}$ and, as stated in the introduction, we will also use the notation $X\simeq_\Q Y$ to mean $X$ and $Y$ have the same rational homotopy type.  That is, $X\simeq_\Q Y$ means there is a zigzag of maps $X\rightarrow Z_1 \leftarrow Z_2\rightarrow ... \leftarrow Y$ each of which induces isomorphisms on all rational homotopy groups.  We note that if $X\simeq_\Q Y$, then $H^\ast(X;\Q)\cong H^\ast(Y;\Q)$.

In general, if $H^\ast(X;\Q)\cong H^\ast(Y;\Q)$, there is no reason to expect that $X\simeq_\Q Y$.  However, there is a certain class of spaces, the so called formal spaces \cite[pg. 156]{FeHaTh}, where this implication does hold.  For example, spheres and projective spaces are formal and the product of two spaces is formal if one of them has rational homology of finite type.  In particular, $S^3\times S^3$, $\mathbb{H}P^m\times S^2$, and $\mathbb{H}P^m\times S^3$ are formal.  Thus, to establish Theorem \ref{rattype}, it is enough to show that if $X$ has $4$-periodic rational cohomology ring which is not singly generated, then the cohomology ring is isomorphic to that of $S^3\times S^3$, $\mathbb{H}P^m\times S^2$ or $\mathbb{H}P^m\times S^3$.

We are now ready to prove Theorem \ref{rattype}, which we recall for convenience.

\mainone*

\begin{proof}\label{mainoneproof}  As mentioned above, it is sufficient to show $H^\ast(X;\Q)$ is isomorphic to the rational cohomology ring of $S^3\times S^3$ or $\mathbb{H}P^m\times S^l$ with $l\in\{2,3\}$.  We note that a rationally elliptic space always has non-negative Euler characteristic and has positive Euler characteristic iff all odd rational Betti numbers vanish \cite[Proposition 32.10]{FeHaTh}.

First, assume $n = 6$ and that the fourth rational Betti number, $b_4$, is equal to $0$.  Then Poincar\'e duality implies $b_2 = 0$ as well, so $0 \leq \chi(X) = 2 - b_3$.  If $\chi(X)>0$, then $b_3 = 0$ so $X$ is rationally $S^6$, so we assume $b_3 = 2$.  Then Poincar\'e duality implies $X$ has the rational cohomology ring of $S^3\times S^3$.

We next assume $b_4 > 0 $, so, by periodicity, $b_4 = 1$.  Poincar\'e duality then implies $b_2 =1$ and thus, $0\leq \chi(X) = 4-b_3$.  If $b_3\leq 3$, $\chi(X) > 0$ so $b_3 = 0$ and $X$ has the rational cohomology groups of $S^4\times S^2$.  If the square of a non-zero element of $H^2(X;\Q)$ is non-zero, it follows that the rational cohomology ring of $X$ is isomorphic to that of $\mathbb{C}P^3$.  Thus, we assume the square of any element in $H^2(X;\Q)$ is $0$.  Using Poincar\'e duality, it is now easy to see the rational cohomology ring is isomorphic to that of $S^4\times S^2 = \mathbb{H}P^1\times S^2$.

Thus, we are left with the case $b_3 = 4$.  However, we now show this cannot occur for a rationally elliptic manifold.  By the rational Hurewicz theorem \cite{KK}, the map $\pi_3(X)_\Q \rightarrow H_3(X;\Q)$ is surjective, so we must have $\dim\pi_3(X)_\Q \geq 4$.  But, according to \cite[pg. 434]{FeHaTh}, for any rationally elliptic space, $\sum_{k} (2k+1)\dim \pi_{2k+1}(X)_\Q\leq 2n-1$.  Since $n = 6$, this is a contradiction, concluding the case of $ n  =6$.

\

We now assume $n > 6$.  Then it is easy to see that if $b_4= 0$ and $X$ is rationally $4$-periodic, then $X$ has the cohomology ring of $S^n$.  Hence, we will assume $b_4 >0 $ so, by periodicity, $b_4=1$.

If the dimension of $X$ is of the form $4m$ or $4m+1$, then Wilking \cite[Proposition 7.13]{Wi2} has proven that $4$-periodicity implies the cohomology ring is singly generated.

Now, suppose that $X$ has dimension $4m+3$.  Then $b_{4m} =1$ by periodicity, so $b_3 = 1$ by Poincar\'e duality.  Periodicity and Poincar\'e duality then imply the rational cohomology ring is isomorphic to that of $\mathbb{H}P^m \times S^3$.

Next, if $X$ has dimension $4m+2\geq 10$, then Amann and Kennard \cite[Lemma 3.1]{AK1} prove $\chi(X) > 0$ so all odd rational Betti numbers vanish.  Using Poincar\'e duality, it follows that $b_{2k} = 1$ for all $k$.  As in the $n=6$ case above, if an element of $H^2(X;\Q)$ has a non-zero square, then the rational cohomology ring of $X$ is isomorphic to that of $\mathbb{C}P^{2m+1}$, so is singly generated.  On the other hand, if every element of $H^2(X;\Q)$ squares to $0$,  it is easy to see the rational cohomology ring is isomorphic to that of $\mathbb{H}P^m\times S^2$.

\end{proof}

\subsection{Biquotients and their classification}

As mentioned in the introduction, a biquotient is any manifold which is diffeomorphic to the quotient of a homogeneous space by a free isometric action.  There is an alternative characterization given in terms of Lie groups.  Suppose $f=(f_1,f_2):U\rightarrow G\times G$ is a homomorphism.  This defines an action of $U$ on $G$ by $u\ast(g) = f_1(u) g f_2(u)^{-1}$.  When this action is free, the orbit space, denoted $G\bq U$ is also called a biquotient.  In the special case when $U = U_1\times U_2$ with each factor embedded into a single factor of $G\times G$, we write $G\bq U \cong U_1\backslash G/U$.  We note that the diffeomorphism type of $G\bq U$ depends only on the conjugacy class of the image of $U$ in $G\times G$ \cite[Proposition 2.2]{DeV2}.

In \cite{To1}, Totaro provides the framework for the classification of biquotients.  By \cite[Lemma 3.3]{To1}, together with the fact that we allow our actions to have finite ineffective kernel, we can and will always assume our biquotients $M = G\bq H$ are \textit{reduced}, meaning $G$ is simply connected, $H\cong H'\times T^k$ is product of a simply connected compact Lie group $H'$ and a torus, and no simple factor of $H$ acts transitively on any simple factor of $G$.

To describe Totaro's main result, first recall that every Lie group is rationally a product of odd spheres, $G\simeq_\Q S^{2d_1 -1} \times ... \times S^{2d_k -1}$ where $k$ is the rank of the group.  The integers $d_i$ are referred to as degrees of $G$.  Table \ref{table:degrees} lists all the simple Lie groups together with their degrees.

\begin{table}[h]

\begin{center}

\begin{tabular}{|c|c|c|}

\hline

Group & Restriction & Degrees \\

\hline

\hline

$SU(n)$ & $n\geq 2$ & $ 2,3,...,n$\\

$SO(2n+1)$ & $n\geq 3$ & $2,4,...,2n$\\

$Sp(n)$ & $n\geq 2$ & $2,4,...,2n$\\

$SO(2n)$ & $n\geq 4$ & $2,4,...,2n-2, n$\\

\hline

$G_2$ & &$2,6$\\

$F_4$ & &$2,6,8,12$ \\

$E_6$ & & $2,5,6,8,9,12$\\

$E_7$ & & $2,6,8,10,12,14,18$\\

$E_8$ & &$2,8,12,14,18,20,24,30$ \\

\hline

\end{tabular}

\caption{Degrees of Lie groups}\label{table:degrees}

\end{center}

\end{table}

Now, suppose $G_1$ is a simple factor of $G$.  We say $G_1$ \textit{contributes degree $d$} to $G\bq H$ if in the long exact sequence of rational homotopy groups associated to the fibration $H\rightarrow G\rightarrow G\bq H$, the induced map $\pi_{2d-1}(H)_\Q\rightarrow \pi_{2d-1}(G_1)_\Q$ is not surjective.

We briefly note that because we allow the $H$ action on $G$ to have finite ineffective kernel $K$, we really only have a fibration of the form $H/K\rightarrow G\rightarrow G\bq H$.  However, since $H$ is a finite cover of $H/K$, the rational homotopy groups of $H$ and $H/K$ are canonically isomorphic, so the above definition makes sense with respect to this slight abuse of notation.

We may now state Totaro's main classification theorem.

\begin{theorem}\label{totthm}\cite[Theorem 4.8]{To1}

Suppose $M = G\bq H$ is a reduced biquotient.  Let $G_1$ denote a simple factor of $G$.  Then one of the following happens.

1.  $G_1$ contributes its highest degree to $M$.

2.  $G_1$ contributes its second highest degree to $M$, and there is a simple factors $H_1$ of $H$ for which $H_1$ acts only on one side of $G_1$, with $G_1/H_1$ isomorphic to $SU(n)/Sp(n)$ with $n\geq 2$,  $Spin(7)/G_2 = S^7$, $Spin(8)/G_2 = S^7\times S^7$, or $E_6/F_4$.  The second highest degrees are given by $2n-1, 4, 4, 9$ respectively.

3.  $G_1\cong Spin(2n)$ and there is a simple factor $H_1$ of $H$ with $H_1 = Spin(2n-1)$, acting only on one side of $G_1$ via the standard inclusion, so $G_1/H_1 = S^{2n-1}$.  In this case, $G_1$ contributes degree $n$ to $M$.

4.  $G_1 = SU(2n+1)$ contributes degrees $2,4,6,...,2n$ to $M$ and there is a simple factor $H_1$ of $H$ with $H_1 \cong SU(2n+1)$ acting on $G_1$ by $h\ast g = hgh^t$.

\end{theorem}

In particular, each simple factor of $G$ contributes at least one degree.  Since are interested in biquotients with the rational homotopy type of $\mathbb{H}P^m\times S^{\ell}$, we therefore may assume $G$ has at most $2$ factors.

\section{\texorpdfstring{Classification when $M\simeq_\Q \mathbb{H}P^m\times S^3$}{Classification when M ~ HPm x S3}}\label{S3}

In this section, we classify biquotients $M = G\bq H$ for which $M\simeq_\Q \mathbb{H}P^m\times S^3$ up to diffeomorphism, proving Theorem \ref{bclass}.  In Subsection \ref{S3action}, we show that every simply connected biquotient which is rationally $\mathbb{H}P^m\times S^3$ is diffeomorphic to one of the spaces listed in Theorem \ref{bclass}.  In Subsection \ref{S3dif}, we show the spaces listed in Theorem \ref{bclass} are distinct up to diffeomorphism.

\subsection{Classification of actions}\label{S3action}

In this section, we assume that $M\simeq_\Q \mathbb{H}P^m \times S^3$.  Hence, $\dim \pi_3(M)_\Q = \dim \pi_4(M)_\Q = \dim \pi_{4m+3}(M)_\Q = 1$  and $\pi_k(M)_\Q = 0$ for all other $k$.

\begin{proposition}\label{classstart}  Suppose $M\simeq_Q \mathbb{H}P^m \times S^3$ and that $M=G\bq H$ is a reduced biquotient.

(1)  $G$ and $H$ have the same number of simple factors.

(2)  The rank of $H$ is one less than the rank of $G$.

\end{proposition}

\begin{proof}

Both assertions follow from the long exact sequence in rational homotopy groups associated to the fibration $H\rightarrow G\rightarrow G\bq H$, together with the well known facts that rank of $\pi_3$ of a Lie group is equal to the number of simple factors and that sum of the dimensions of the odd degree rational homotopy groups is equal to the rank.
\end{proof}

We have already argued that $G$ has at most two simple factors.  We now show that $G$ can not be simple.

\begin{proposition}\label{simplecase}  Suppose $M=G\bq H\simeq_\Q \mathbb{H}P^m\times S^3$ is a reduced biquotient.  Then $G$ cannot be simple.
\end{proposition}

\begin{proof}

Assume $G$ is simple.  By Proposition \ref{classstart}(1), $H$ is simple as well, so $\pi_3(G)_\Q\cong \pi_3(H)_\Q \cong \Q$.  Now, consideration of the long exact sequence in rational homotopy groups shows that $\pi_{2i-1}(H)_\Q \cong \pi_{2i-1}(G)_\Q$ with $i\neq 2m+2$ and $i\geq 3$.

On the other hand, when $i = 2m+2$, $\pi_{4m+3}(M)_\Q \cong\Q$, so $\pi_{4m+3}(G)_\Q\cong \Q$ but $\pi_{4m+3}(H)_\Q = 0$.

As $G$ is rationally a product of spheres, we see $H^\ast(G;\Q) \cong H^\ast(H\times S^{4m+3};\Q)$.  But the pairs $(G,H)$ for which this happens are classified in \cite[Section 2]{KZ} and, in particular, for any such pair $(G,H)$, any biquotient $G\bq H$ has singly generated rational cohomology.

\end{proof}

We may now assume $G = G_1\times G_2$ and $H = H_1\times H_2$.  From Theorem \ref{totthm}, we may assume $G_1$ contributes degree $2m+2$ and $G_2$ contributes degree $2$.  If $2$ is the highest degree of $G_2$, then $G_2 = SU(2)$.  If $2$ is not the highest degree of $G_2$, then, from Theorem \ref{totthm}, we see that $G_2$ must come from case (4), so $G_2 = H_2 = SU(3)$ with $H_2$ acting on $G_2$ by $A\ast B = ABA^t$.  We now show the case where $G_2 = H_2 = SU(3)$ cannot occur.

\begin{proposition}

If $G = G_1\times SU(3)$, $H = H_1\times SU(3)$ with the $SU(3)$ factor of $H$ acting on the $SU(3)$ factor of $G$ via $A\ast B = ABA^t$, then $G\bq H$ does not have the rational homotopy type of $\mathbb{H}P^m\times S^3$.

\end{proposition}

\begin{proof}We proceed by contradiction.  The homomorphism $SU(3)\rightarrow SU(3)^2$ given by $A\mapsto(A,\overline{A})$ defines the action $A\ast B = ABA^t$.  This homomorphism has image a maximal connected subgroup and thus, $H_1$ must act effectively freely on the first factor $G_1$.  Since $H$ has rank $1$ less than that of $G$, it follows that $H_1$ also has rank one less than that of $G_1$.  In particular, $G_1\bq H_1$ is odd dimensional.  Further, just as in the proof of Proposition \ref{simplecase}, the degrees of $H_1$ and $G_1$ must agree, except for a single degree of $G_1$.  It follows from \cite[Section 2]{KZ} that $G_1\bq H_1$ is a simply connected rational sphere.

Now, consider the restriction of the $H$ action to $H_1\times SO(3)\subseteq H_1\times SU(3)$.  The projection of this action to the $SU(3)$ factor of $G$ fixes the identity, and so $H_1\times SO(3)$ must act effectively freely on $G_1$.  In particular, $SO(3)$ acts effectively freely on $G_1\bq H_1$.  In fact, since $SO(3)$ has no non-trivial normal subgroup, it must act freely on $G_1\bq H_1$, so we have a principal $SO(3)$ bundle \begin{equation} SO(3)\rightarrow G_1\bq H_1\rightarrow G_1\bq (H_1\times SO(3)) \tag{$\ast$}.\end{equation}  From the long exact sequence in rational homotopy groups together with the rational Gysin sequence, we see that $G_1\bq(H_1\times SO(3))\simeq_\Q \mathbb{H}P^m$ with $\pi_2$ non-trivial.

Now, biquotients with singly generated rational cohomology rings were classified in \cite{KZ}: the only example of a rational $\mathbb{H}P^m$ with non-trivial $\pi_2$ is $\mathbf{G}_2/SO(4)$.  It follows that $G = \mathbf{G}_2\times SU(3)$ and $H = SU(2)\times SU(3)$, up to cover.

But Eschenburg \cite[pg. 122]{Es2} classified biquotient actions of maximal rank on simple Lie groups and, in particular, showed that $\mathbf{G}_2$ admits no free two-sided action.  Hence, $H_1 = SU(2)$ and $SO(3)\subseteq SU(3)$ must both act on the same side of $\mathbf{G}_2$.  But, the only non-trivial homomorphism $SU(3)\rightarrow \mathbf{G}_2$ has maximal image, so $H_1$ cannot act on $\mathbf{G}_2$.  This contradiction shows the case $G = G_1\times SU(3) $ cannot occur.

\end{proof}

We henceforth assume $G = G_1\times SU(2)$ and $H= H_1\times H_2$.  Because $H$ does not contain a torus factor, \cite[Proposition 4.1]{DeV2} shows the projection of the $H$ action onto $G_1$ must be effectively free.  Further, the projection onto the $SU(2)$ factor can only be non-trivial if one factor of $H$ is isomorphic to $SU(2)$, which then acts by conjugation.

The biquotient $G_1\bq H$ is simply connected because $G_1$ is simply connected and $H$ is connected.  Further, it is rationally $2$-connected since $\pi_2(G_1) = 0$ and $\pi_1(H)$ is finite.  Moreover, it is rationally $3$-connected.  To see this note that $H$ has the same rank as $G_1$, so $G_1\bq H$ has positive Euler characteristic \cite[Theorem 5.1]{Si1}.  But, as mentioned in the proof of Theorem \ref{rattype}, for rationally elliptic topological spaces, the Euler characteristic is positive iff all odd Betti numbers vanish.

Because the $H$ action on $G_1$ is effectively free, we have a bundle $H/K\rightarrow G_1\rightarrow G_1\bq H$, where $K$ is the ineffective kernel of the action.  The $H$ action on the $S^3 = SU(2)$ factor of $G$ gives an associated bundle $S^3\rightarrow G_1\times_H S^3 = G\bq H\rightarrow G_1\bq H$.  We use the associated bundle to study the topology of $G_1\bq H$.

\begin{proposition}\label{case3}  $G_1\bq H$ has the rational cohomology ring of $\mathbb{H}P^m$.
\end{proposition}

\begin{proof}  Consider the rational Gysin sequence associated to the bundle $S^3\rightarrow G\bq H\rightarrow G_1\bq H$ and recall that $G\bq H$ is rationally $\mathbb{H}P^m\times S^3$.  Since $H^3(G_1\bq H;\Q) = 0$, the rational Gysin sequence shows that $H^3(G\bq H)\cong \Q$ must inject into $H^0(G_1\bq H;\Q) \cong \Q$, which implies the Euler class of the bundle is $0$.  It now follows that the projection map $G\bq H\rightarrow G_1\bq H$ embeds $H^\ast(G_1\bq H;\Q)$ as a subalgebra of $H^\ast(G\bq H;\Q)$.  This shows that a generator of $H^4(G_1\bq H;\Q)$ has non-trivial powers until $x^{m+1} = 0$.  Thus, $H^\ast(G_1\bq H;\Q)\cong H^\ast(\mathbb{H}P^m;\Q)$.

\end{proof}

In \cite{KZ}, Kapovitch and Ziller classify biquotients with the rational cohomology ring of $\mathbb{H}P^m$: such a biquoient is diffeomorphic to $\mathbb{H}P^m$, $\mathbf{G}_2/SO(4)$, or $N^m:= \Delta SU(2)\backslash SO(2m+3)/SO(2m+1)$, where the notation $\Delta SU(2)$ refers to the composition $SU(2)\rightarrow SO(4)\rightarrow SO(2m+3)$ with $SO(4)\rightarrow SO(2m+3)$ the block diagonal embedding $A\mapsto \diag(A, A,..., A, 1)$.  The examples are distinguished up to homotopy by their cohomology rings.

When $G_1\bq H$ is not diffeomorphic to $\mathbb{H}P^m$, we see that, up to cover, $H = H_1\times SU(2)$.  One the other hand, a reduced biquotient which is diffeomorphic to $\mathbb{H}P^m$ must also have $H = H_1\times SU(2)$, as follows from \cite[pg. 264]{Oni} in the homogeneous case and \cite[Table 101]{Es2} in the inhomogeneous case.  Hence, we may assume $H_2 = SU(2)$.

The following theorem summarizes the structure of $G\bq H\sim_Q \mathbb{H}P^m\times S^3$.

\begin{theorem}\label{hpms3}

Every compact simply connected biquotient $M = G\bq H$ with $M\simeq_\Q \mathbb{H}P^m\times S^3$ is diffeomorphic to a biquotient having the form $$G\bq H = (G_1\times SU(2))\bq (H_1\times SU(2))$$ where the projection of the $H$ action onto $G_1$ is effectively free with quotient a rational $\mathbb{H}P^m$.  In addition, either the projection of the $H$ action on the $SU(2)$ factor of $G$ is trivial or at most one factor of $H$, isomorphic to $SU(2)$, acts by conjugation.

\end{theorem}

If $m=1$, $G\bq H$ is $7$-dimensional.  Such biquotients were classified in \cite[Propositions 4.27 and 4.28]{DeV2}.  In particular among those with the rational homotopy of $\mathbb{H}P^1\times S^3$, the homotopy and diffeomorphism classification agree and there are precisely three diffeomorphism types.

We now assume $m\geq 2$.  With this restriction, the only case where $H_1 = SU(2)$ occurs when $G_1 = \mathbf{G}_2$.  Thus, in this case, we get three biquotients:  the trivial product $(\mathbf{G}_2/SO(4)) \times S^3$, as well as two biquotients $(\mathbf{G}_2/H_1)\times_{H_2} S^3$ and $(\mathbf{G}_2/H_2)\times_{H_1} S^3$.  In  the remaining cases, we end up with four forms of biquotients:  the trivial products $\mathbb{H}P^m\times S^3$ and $N^m\times S^3$, as well as the non-trivial bundles $S^{4m+3}\times_{SU(2)} S^3$ and $(SO(2m+3)/SO(2m+1))\times_{SU(2)}S^3$.

We have now shown half of Theorem \ref{bclass}, that each such biquotient is diffeomorphic to one of the listed spaces.  In the next section, we prove these biquotients are distinct up to diffeomorphism.

\subsection{Diffeomorphism classification}\label{S3dif}

In this section, we show each of the biquotients listed in Theorem \ref{bclass} is distinct up to diffeomorphism.

Recall that, in the proof of Proposition \ref{case3}, we showed that if $M = G\bq H\simeq_\Q \mathbb{H}P^m\times S^3$, then we have a bundle $S^3\rightarrow G\bq H \rightarrow G_1\bq H$ with $G_1\bq H$ an element of $\{\mathbb{H}P^m, N^m, \mathbf{G}_2/SO(4)\}$ and with Euler class $0$.  Thus, the Gysin sequence breaks into short exact sequences of the form $$0\rightarrow H^\ast(G_1\bq H)\rightarrow H^\ast(G\bq H) \rightarrow H^{\ast - 3}(G_1\bq H)\rightarrow 0.$$

If $G_1\bq H\neq \mathbf{G}_2/SO(4)$, then for any $\ast$, either $H^\ast(G_1\bq H) = 0$ or $H^{\ast - 3}(G_1\bq H) = 0$.  It follows that the integral cohomology ring $H^\ast(G\bq H)$ is isomorphic that of $(G_1\bq H)\times S^3$.  If, on the other hand, $G_1\bq H = \mathbf{G}_2/SO(4)$, then we see the projection map induces an isomorphism $\mathbb{Z}/2\mathbb{Z} = H^3(G_1\bq H)\rightarrow H^3(G\bq H)$, so $H^\ast(G\bq H)$ has torsion. Since the cohomology rings of the three rational $\mathbb{H}P^m$s distinguish them, we have shown the following proposition.

\begin{proposition}\label{s3hom}

Suppose $$B_1,B_2\in\{\mathbb{H}P^m, N^m, \mathbf{G}_2/SO(4)\}$$ with $B_1\neq B_2$  Let $M_i$ be the total space of a linear $S^3$ bundle over $B_i$ with Euler class $0$.  Then $M_1$ is not homotopy equivalent to $M_2$
\end{proposition}

Thus, in terms of the diffeomorphism classification, each of cases (1), (2), and (3) of Theorem \ref{bclass} can be treated separately.

We recall that when $B = \mathbb{H}P^m$ or $N^m$, the bundle $S^3\rightarrow G\bq H\rightarrow B$ is associated to a principal bundle  $$SU(2)\rightarrow G_1\bq H_1 \stackrel{\pi}{\rightarrow} G_1\bq H.$$

When $B = \mathbf{G}_2/SO(4)$, there are two relevant principal bundles.  The two normal $SU(2)$s in $SO(4)$ have Dynkin indices $1$ and $3$ \cite{KZ}, so we denote them $SU(2)_1$ and $SU(2)_3$.  Then $SU(2)_3$ acts effectively free on $\mathbf{G}_2/SU(2)_1$ with ineffective kernel $\{\pm I\}\in SU(2)$ giving a principal $SO(3)$ bundle $\mathbf{G}_2/SU(2)_1\rightarrow \mathbf{G}_2/SO(4)$.  The same argument also gives a principal $SO(3)$-bundle $\mathbf{G}_2/SU(2)_3\rightarrow \mathbf{G}_2/SO(4).$

We let $F\in \{SU(2), SO(3)\}$ denote the relevant fiber of the principal bundle.  Now, for any choice of $B$, the action of $F$ on $S^3$ is either trivial or by conjugation.  Hence, we must distinguish $B\times S^3$, when the action is trivial, from $(G_1\bq H_1)\times_{F} S^3$ where the $F$ action on $S^3$ is by conjugation.  We do this by computing Pontrjagin classes.

Our main tool for computing the first Pontrjagin class of a biquotient of the form $(G_1\bq H_1)\times_F S^3 $ is the following proposition, which can be found, for example, in \cite[pg. 383]{Ge}.

\begin{proposition}\label{totbun}  Suppose $S^n\rightarrow E\xrightarrow{\pi} B$ is an oriented linear sphere bundle and let $\xi = \mathbb{R}^{n+1}\rightarrow \tilde{E}\xrightarrow{\tilde{\pi}} B$ denote the corresponding vector bundle.  Then $T\tilde{E} \cong \tilde{\pi}^\ast(\xi)\oplus \tilde{\pi}^\ast(TB)$.  In addition, the inclusion $i:E\rightarrow \tilde{E}$ gives $i^\ast(T\tilde{E}) \cong TE \oplus 1$, with $1$ denoting the rank $1$ trivial bundle.  In particular, $$TE \oplus 1 \cong \pi^\ast(\xi)\oplus \pi^\ast(TB).$$
\end{proposition}

Using this, we now show $B\times S^3$ and $(G_1\bq H_1)\times_F S^3$ have different diffeomorphism types.

\begin{proposition}\label{s3diff}For any $$B\in \{\mathbb{H}P^m, N^m, \mathbf{G}_2/SO(4)\},$$ the trivial bundle $B \times S^3$ and the non-trivial bundle $(G_1\bq H_1)\times_F S^3$ have distinct Pontrjagin classes.

\end{proposition}

\begin{proof}

The $F$ action on $S^3$ is linear, so extends to a representation on $\mathbb{R}^4$.  Let $\xi$ denote the rank 4 vector bundle $(G_1\bq H_1)\times_F \mathbb{R}^4\rightarrow G_1\bq H$.

From Proposition \ref{totbun}, we see $$p_1((G_1\bq H_1)\times_{SU(2)} S^3) =  \pi^\ast(p_1(\xi)) + \pi^\ast(p_1(TB)),$$ and so we need only show $\pi^\ast(p_1(\xi))\neq 0$.  Since the associated bundle $S^3\rightarrow G\bq  H\rightarrow B$ has trivial Euler class, $\pi^\ast$ is injective on $H^4$, so it is enough to show $p_1(\xi)\neq 0$.

\

Consider the universal principal $F$ bundle $EF\rightarrow BF$.  Then we have a commutative diagram of fibrations \begin{diagram}F & \rTo & G_1\bq H_1 & \rTo & G_1\bq H \\ & & \dTo & & \dTo^{\phi} \\ F & \rTo & EF & \rTo & BF \end{diagram} where $\phi:G_1\bq H\rightarrow BF$ is the classifying map.  This gives rise to a map of associated bundles \begin{diagram} \mathbb{R}^4 & \rTo & (G_1\bq H_1)\times_F \mathbb{R}^4 & \rTo & G_1\bq H \\ & & \dTo & & \dTo^{\phi}\\ \mathbb{R}^4 & \rTo & EF\times_{F}\mathbb{R}^4 & \rTo & BF.\end{diagram}  Thus, if $\eta$ denotes the vector bundle $\mathbb{R}^4\rightarrow EF\times_{F} \mathbb{R}^4\rightarrow BF$, then $\xi$ is the pull back of $\eta$, $\xi = \phi^\ast \eta$.  Hence, in order to show $p_1(\xi)\neq 0$, it is enough to show $p_1(\eta)\neq 0$ and that $\phi^\ast$ is non-zero on $H^4$.

\

Let $i:S^1\rightarrow F$ denote the inclusion of a maximal torus.  Then, for either choice of $F$, it is easy to see that the induced map $i^\ast:H^4(BF)\rightarrow H^4(BS^1)$ is an isomorphism.

It is well known \cite[Theorem 10.3]{BH1} that the associated vector bundle $\eta$ over $BF$ has total Pontrjagin classes given by $i^\ast \prod (1+\beta^2)$ where the product is over all weights of the representation.  Since the conjugation action has non-trivial weights, it follows immediately that $p_1(\eta) = \sum \beta^2 \neq 0$.

Thus, to complete the proof, we need only show $\phi^\ast:H^4(BF)\rightarrow H^4(G_1\bq H)$ is non-zero.  If it is the $0$ map, naturality of Serre spectral sequences implies that all differentials in the spectral sequence of the principal $F$ bundle $G_1\bq H_1\rightarrow G_1\bq H$ vanish, which, in turn, implies that $G_1\bq H_1$ has the rational cohomology ring of $\mathbb{H}P^m\times S^3$.  Since $G_1$ is simple, this contradicts Proposition \ref{simplecase}.

\end{proof}

Together, Propositions \ref{s3hom} and \ref{s3diff} almost completely classify biquotients which are rationally $\mathbb{H}P^m\times S^3$ up to diffeomorphism.  The only remaining task is to distinguish the two nontrivial bundles over $\mathbf{G}_2/SO(4)$.  Thus, once we prove the following proposition, we will have completed the proof of Theorem \ref{bclass}.

\begin{proposition}\label{g2class}  The first Pontrjagin classes of the biquotients $$(\mathbf{G}_2/SU(2)_1)\times_{SO(3)} S^3\text{ and }(\mathbf{G}_2/SU(2)_3)\times_{SO(3)} S^3$$ are different.

\end{proposition}

\begin{proof}

We first recall the torsion in the cohomology ring of $G_2/SO(4)$ is $2$-torsion and that $H^4(G_2/SO(4))\cong \Z$ \cite[pg. 242]{IsTo}.

We have two commutative diagrams \begin{diagram}SO(3) & \rTo & \mathbf{G}_2/SU(2)_i & \rTo & \mathbf{G}_2/SO(4)\\ & & \dTo & & \dTo^{\phi_i} \\ SO(3)& \rTo & ESO(3) & \rTo & BSO(3) \end{diagram} for $i =1, 3$.  Following the proof of Proposition \ref{s3diff}, we need only show $\phi_i^\ast:\mathbb{Z} = H^4(BSO(3))\rightarrow H^4(\mathbf{G}_2/SO(4))$ for $i =1,3$ are different maps.

When $i = 1$, $\mathbf{G}_2/SU(2)_1$ is diffeomorphic to the unit tangent bundle of $S^6$, which is $5$-connected.  In particular, the differential $d:\Z\cong H^3(SO(3);\Z)\rightarrow H^4(\mathbf{G}_2/SO(4))$ must be an isomorphism.  From naturality of the Serre spectral sequence, the map $H^4(BSO(3))\rightarrow H^4(G_2/SO(4))$ must be surjective, so is multiplication by $\pm 1$.

On the other hand, when $i = 3$, $\mathbf{G}_2/SU(3)_3$ is $2$-connected with $\pi_3(\mathbf{G})_2/SU(2)_3 \cong \Z/3\Z$, and hence $H^4(\mathbf{G}_2/SU(3)_3)$ contains $3$-torsion.  Since the only torsion in the cohomology rings of $SO(3)$ and $\mathbf{G}_2/SO(4)$ are $2$-torsion, it follows from naturality of the Serre spectral sequence that the map $H^4(BSO(3))\rightarrow H^4(\mathbf{G}_2/SO(4))$ is multiplication by $\pm 3$.
\end{proof}

\section{\texorpdfstring{Classification when $M\simeq_\Q \mathbb{H}P^m\times S^2$}{Classification when M ~ HPm x S2}}\label{S2}

We now study biquotients $M = G\bq H$ with $M\simeq_\Q \mathbb{H}P^m\times S^2$.  We begin by studying the possible actions of $H$ on $G$, finding a structure theorem analogous to Theorem \ref{hpms3}.  Then, we partially classify the resulting diffeomorphism types.

\subsection{Classification of actions}

We prove the following structure theorem by reducing classification of biquotients $M$ which are rationally $\mathbb{H}P^m\times S^2$ to those which are rationally $\mathbb{H}P^m\times S^3$.

\begin{theorem}\label{hpms2}

Every compact simply connected biquotient $M = G\bq H$ with $M\simeq_\Q \mathbb{H}P^m\times S^2$ is diffeomorphic to a biquotient having the form $G\bq H = (G_1\times SU(2))\bq (H_1\times SU(2)\times S^1)$.  Further, the projection of the $H_1\times SU(2)$ action onto $G_1$ is effectively free with quotient a rational $\mathbb{H}P^m$ and, in addition, the projection of the $H_1\times SU(2)$ action on the $SU(2)$ factor of $G$ is trivial.  Finally, the circle factor of $H$ acts, up to ineffective kernel, as the Hopf map on the $SU(2)$ factor of $G$.

\end{theorem}

\begin{proof}

Suppose $M = G\bq H$ is a reduced biquotient with $M\simeq_\Q \mathbb{H}P^m\times S^2$.  Since $\pi_2(M)_\Q\cong \Q$, it follows that $H = H'\times S^1$ with $H'$ semi-simple.  This gives a principal $S^1$-bundle $G\bq H'\rightarrow G\bq H$.  It follows from the long exact sequence in rational homotopy groups that $G\bq H' \simeq_\Q \mathbb{H}P^m\times S^3$.  Applying Theorem \ref{hpms3}, $G = G_1\times SU(2)$, $H' = H_1\times SU(2)$, where the projection of the $H'$ action to $G_1$ is effectively free with quotient a rational $\mathbb{H}P^m$.

We claim that the projection of the $H'$ action to the $SU(2)$ factor of $G$ is trivial.  Indeed, from \cite[Proposition 4.1]{DeV2}, if the $H'$ action on $SU(2)$ is non-trivial, then the projection of the $H = H'\times S^1$ action to $G_1$ must be effectively free.  But there is no free $S^1$ action on a rational $\mathbb{H}P^m$: every circle action has fixed points since $\chi(\mathbb{H}P^m) = m+1 \neq 0$.

Finally, since the diagonal action of $S^1$ on $(G_1\bq H')\times SU(2)$ must be free, the projection of the $S^1$ action to the $SU(2)$ factor of $G$ must be effectively free, so it is, up to ineffective kernel, the Hopf action.

\end{proof}

In fact, under the hypothesis of Theorem \ref{hpms2}, $M$ is diffeomorphic to a space of the form $(G_1\bq (H_1\times SU(2)))\times_{S^1} S^3$ where the action on $S^3$ is effective.  For, if $z\in S^1$ acts trivially on $S^3$, then it fixes a point of $(G_1\bq (H_1\times SU(2))) \times S^3$ because the $S^1$ action on $G_1\bq (H_1\times SU(2))$ has a fixed point.  Because the action is effectively free, we conclude that $z$ fixes all points of $(G_1\bq (H_1\times SU(2)))\times S^3$.  Thus, dividing $S^1$ by its ineffective kernel, we obtain an effective action of $S^1$ on $(G_1\bq (H_1\times SU(2)))\times S^3$ as desired.

Because $M$ is diffeomorphic to a space of the form $(G_1\bq (H_1\times SU(2)))\times_{S^1} S^3$, we see that $M$ has the structure of the total space of a bundle over $S^2$ with fiber $\mathbb{H}P^m, N^m,$ or $\mathbf{G}_2/SO(4)$.  In fact, the bundle structure is associated to the Hopf bundle $S^1\rightarrow S^3\rightarrow S^2$ via the $S^1$ action on $G_1\bq(H_1\times SU(2))$.  This is in contrast to the case where $M\simeq_\Q \mathbb{H}P^m\times S^3$, where $M$ naturally has the structure of a bundle with fiber $S^3$ and base one of $\mathbb{H}P^m, N^m,$ or $\mathbf{G}_2/SO(4)$.

\subsection{Partial diffeomorphism classification}\label{partdif}

In this section, we investigate the topology of biquotients $M=G\bq H$ with $M\simeq_\Q \mathbb{H}P^m\times S^2$, proving Theorem \ref{cclass}.  We have already shown that such an $M$ is a bundle over $S^2$ with fiber over $B$, $$B\in\{\mathbb{H}P^m, N^m, \mathbf{G}_2/SO(4)\}.$$  When $m = 1$, such biquotients are classified in \cite[Section 4.1]{DeV2}.  In particular, the homotopy and diffeomorphism classifications coincide, and there are precisely two diffeomorphism types.  Thus, we now assume $m\geq 2$.

The differentials in the Serre spectral sequence associated to $B\rightarrow M\rightarrow S^2$ all vanish for trivial reasons and there are no extension problems.  If $B = \mathbb{H}P^m$ or $N^m$, it follows that $H^\ast(M)$ is torsion free.  However, the inclusion $B\rightarrow M$ induces ring isomorphisms $H^{4i}(M)\rightarrow H^{4i}(B)$ for $0\leq i\leq m$.  Thus, $H^\ast(M)\cong H^\ast(B\times S^2)$ when $B = \mathbb{H}P^m$ or $N^m$.  In particular, if $M_1$ is a biquotient with $B_1 = \mathbb{H}P^m$ and $M_2$ is a biquotient with $B_2 = N^m$, then $H^\ast(M_1)\not\cong H^\ast(M_2)$.

In addition, when $B = \mathbf{G}_2/SO(4)$, it follows that $H^\ast(M)$ has torsion.  Thus, we have the following analogue of Proposition \ref{s3hom}.

\begin{proposition}
Suppose $$B_1,B_2\in\{\mathbb{H}P^m, N^m, \mathbf{G}_2/SO(4)\}$$ with $B_1\neq B_2$  Let $M_i$ denote a biquotient of the form $B_i \times_{S^1} S^3$ where the action of $S^1$ on $S^3$ is the Hopf map.  Then $M_1$ is not homotopy equivalent to $M_2$
\end{proposition}

Hence, in terms of proving Theorem \ref{cclass}, we may work with one choice of $B$ at a time.  We begin with the case $B = \mathbf{G}_2/SO(4)$.

\begin{proposition}\label{g2diff}

Suppose $M=G\bq H$ is a simply connected biquotient diffeomorphic to $(\mathbf{G}_2/SO(4))\times_{S^1} S^3$ where the $S^1$ action on $S^3$ is the Hopf action.  Then $M$ is diffeomorphic to the product $(\mathbf{G}_2/SO(4))\times S^2$.

\end{proposition}

\begin{proof}

We have already shown $M$ has the form $B\times_{S^1} S^3$ and is an associated bundle to the Hopf bundle.  Equipping $\mathbf{G}_2$ with a bi-invariant metric, the $S^1$ action on $B$ is isometric, so we have a map $S^1\rightarrow Iso_0(B)$, the identity component of the isometry group.  This allows us to extend the structure group of $M$, a $\mathbf{G}_2$ bundle over $S^2$ to $Iso_0(B)$.

But $\mathbf{G}_2/SO(4)$ with bi-invariant metric is a symmetric space, so $Iso_0(B)=\mathbf{G}_2/(Z(\mathbf{G}_2)\cap SO(4)) = \mathbf{G}_2$ since the center of $\mathbf{G}_2$ is trivial.  Now, principal $\mathbf{G}_2$ bundles over $S^2$ are in bijective correspondence with $[S^2,B\mathbf{G}_2]$.  But $\mathbf{G}_2$ is $2$-connected, so $B\mathbf{G}_2$ is $3$-connected.  In particular, the only $\mathbf{G}_2$-principal bundle over $S^2$ is the trivial bundle.

Consider the bundle over $S^2$ given by $\mathbf{G}_2\times_{S^1} S^3$ which is associated to the Hopf bundle.  Left multiplication by $\mathbf{G}_2$ is well defined so this is a principal $\mathbf{G}_2$-bundle, which must therefore be equivariantly trivial.

Then we have \begin{align*} M &\cong B\times_{S^1} S^3\\ &\cong B\times_{\mathbf{G}_2}(\mathbf{G}_2\times_{S^1} S^3)\\ &\cong B\times_{\mathbf{G}_2}(\mathbf{G}_2\times S^2)\\ &\cong B\times S^2.\end{align*}

\end{proof}

Having handled the case $B = \mathbf{G}_2/SO(4)$, we now turn to the case where $B = \mathbb{H}P^m$.  From the classifications of homogeneous spaces \cite[pg. 264]{Oni} and biquotients \cite[Table101]{Es2} which are diffeomorphic to $\mathbb{H}P^m$, a bi-invariant metric on $G_1$ induces, up to scale, the Fubini-Study metric on $G_1\bq (H_1\times SU(2) = \mathbb{H}P^m$.  Thus, the identity component of the isometry group $Iso_0(B)$ is $Sp(m+1)/Z(Sp(m+1)) = Sp(m+1)/(\mathbb{Z}/2\mathbb{Z})$.  Now, because $\pi_1(Iso_0(B)) = \mathbb{Z}/2\mathbb{Z}$, it follows that $[S^2, BIso(B)]$ contains two elements.

Thus, following the proof of Proposition \ref{g2class}, we conclude that for each $m$, there are at most two diffeomorphism types of biquotients of the form $\mathbb{H}P^m\times_{S^1} S^3$, depending on how $S^1$ acts on $\mathbb{H}P^m$.  Of course, the trivial action gives the biquotient $\mathbb{H}P^m\times S^2$.  However, for each $m$, we construct an action of $S^1$ on $\mathbb{H}P^m$ for which the biquotient $\mathbb{H}P^m\times_{S^1}S^3$ is not even homotopy equivalent to $\mathbb{H}P^m\times S^2$:  the two biquotients are distinguished by their Stiefel-Whitney classes, a homotopy invariant \cite{Wu}.

More precisely, for a particular action of $S^1$ on $\mathbb{H}P^m$, we will show that either $w_2(\mathbb{H}P^m\times_{S^1} S^3)\neq 0$ or $w_6(\mathbb{H}P^m\times_{S^1} S^3)\neq 0$.  On the other hand, since $S^2$ is stably parallelizable, $\mathbb{H}P^m\times S^2$ has non-trivial Stiefel-Whitney classes only in dimensions divisible by $4$.

We let $S^1 = \{e^{i\theta}: \theta \in [0,2\pi)\}$ act on $\mathbb{H}P^m = \{[q_0:...:q_m]| q_i \in \mathbb{H}\}$ (where $[q_0:...q_m] \simeq [p_0:...:p_m]$ if there is a $q\in \mathbb{H}$ with $q_j q = p_j$ for all $j$) as follows:  \begin{equation*} e^{i\theta} \ast [q_0:...:q_m] = [e^{i\theta/2} q_0: ... : e^{i\theta/2} q_m]. \tag{$\ast$}.\end{equation*}  We let $C_m = \mathbb{H}P^m\times_{S^1} S^3$ using this action of $S^1$ on $\mathbb{H}P^m$.  We set $H^\ast(C_m;\mathbb{Z}/2\mathbb{Z}) = \mathbb{Z}/2\mathbb{Z}[u_2, u_4]/\langle u_2^2, u_4^{m+1}\rangle$.

\begin{proposition}\label{hpclass}
For each $m$, either $w_2(TC_m)\neq 0$ or $w_6(TC_m)\neq 0$.

\end{proposition}

\begin{proof} We first note that there is a commutative diagram of fibrations: \begin{diagram} \mathbb{H}P^m & \rTo^{j} & C_m &\rTo & S^2\\ \dTo & & \dTo^{i} & & \dTo \\ \mathbb{H}P^{m+1} & \rTo^{j} & C_{m+1} & \rTo & S^2 \end{diagram}

To see this, simply note that $\mathbb{H}P^m$ embeds into $\mathbb{H}P^{m+1}$ as the subset with last coordinate $0$.  Then the $S^1$ action clearly preserves this set and acts as in $(\ast)$.  We point out that on cohomology, $i^\ast$ is an isomorphism except that it is trivial on $H^{4m+6}$.

We recall that $j^\ast:H^\ast(C_m)\rightarrow H^\ast(\mathbb{H}P^m)$ is an isomorphism when $\ast = 4k$.  Further, the normal bundle of $\mathbb{H}P^m\subseteq C_m$ is trivial: a neighborhood of $\mathbb{H}P^m$ takes the form $\mathbb{H}P^m\times \mathbb{R}^2$, coming from a trivializing neighborhood in $S^2$.  It follows that $w_4(C_m)=0$ iff $w_4(\mathbb{H}P^m) = 0$.  But the total Stiefel-Whitney class of $\mathbb{H}P^m$ is given by $(1+j^\ast u_4)^{m+1}$ \cite[Section 15.7]{BH1}.  Thus, $w_4(\mathbb{H}P^m)= w_4(C_m)=0 $ iff $m$ is odd.

\

We next claim that $w_2(TC_1) = u_2$ is non-trivial.  This follows because $\mathbb{H}P^1 = S^4$, and the $S^1$ action on $S^4$ comes from an embedding of $S^1$ into $SO(5) = Sp(2)/Z(Sp(2))$.  This embedding is homotopically non-trivial because the lift to $Sp(2)$ is not a loop.  Then according to \cite[Lemma 8.2.5]{Ge}, $w_2\neq 0$ in this case.  Further, $w_6(C_1) = 0$ since is the reduction mod $2$ of the Euler characteristic, which is $4$.

\

Now, inductively, we have $C_1\subseteq C_2\subseteq ... \subseteq C_m$.  If $i:C_1\rightarrow C_m$ denotes the inclusion, then $i^\ast TC_m \cong TC_1 \oplus \nu_{1,m}$ where $\nu_{i,j}$ is the normal bundle of $C_i$ in $C_j$.  Because of the chain of inclusions, the normal bundle $\nu_{1,m}$ decomposes into a sum of rank $4$ bundles and clearly these rank $4$ bundles, when pulled back to $C_1$, are all isomorphic.  Thus, $i^\ast TC_m \cong TC_1\oplus (m-1)\nu_{1,2} =  TC_1 \oplus \underbrace{\nu_{1,2}\oplus \ldots \oplus \nu_{1,2}}_{m-1}.$

We claim that $w_4(\nu_{1,2}) = u_4$.  Indeed, since the product of any two elements of $H^2(C_m)$ vanishes, the Whitney sum formula gives $w_4(i^\ast TC^2) = w_4(TC^1) + w_4(\nu_{1,2})$.  Since $w_4(TC^m)$ depends on the parity of $m$, $w_4(\nu_{1,2}) = u_4$.

\

We now break into cases depending on whether $w_2(\nu_{1,2}) = 0$ or $w_2(\nu_{1,2}) = u_2$.  First, if $w_2(\nu_{1,2}) = 0$, then $i^\ast w_2(TC_m) = w_2(TC_1)\neq 0$, so $w_2(TC_m)\neq 0$.  In particular it follows in this case that $C_m$ is not homotopy equivalent to $\mathbb{H}P^m\times S^2$.

So, we may assume $w_2(\nu_{1,2}) = u_2 \neq 0$, then \begin{align*} i^\ast w(TC_m) &= w(TC_1)w((m-1)\nu_{1,2})\\ &= (1+u_2)(1+u_2 + u_4)^{m-1}.\end{align*}  The degree two part of $(1+u_2)(1+u_2+u_4)^{m-1}$ is $m u_2$, which is non-trivial if $m$ is odd.  On the other hand, the degree six part of $(1+u_2)(1+u_2 + u_4)^{m-1}$ is $(m-1 + (m-1)(m-2)) u_2 u_4$ which is non-trivial if $m$ is even.  Thus, if $w_2(\nu_{1,2})\neq 0$, either $w_2(TC_m)\neq 0$ or $w_6(TC_m)\neq 0$, so $C_m$ is not homotopy equivalent to $\mathbb{H}P^m\times S^2$.

\end{proof}

We now turn attention to the final case $B = N^m$.  We are, unfortunately, unable to obtain a full classification.  We begin with a theorem of Sullivan \cite[Theorem 13.1]{Su1} (see also \cite[Proposition 2.3(d)]{WZ}).

\begin{theorem}(Sullivan)\label{finitediff} Any class of simply connected manifolds of dimension at least five with isomorphic cohomology rings, the same Pontrjagin classes, and whose minimal models are formal contains at most finitely many diffeomorphism types.

\end{theorem}

From the discussion leading up to the proof of Theorem \ref{rattype}, we have already noted that rationally $4$-periodic biquotients are formal.  Further, at the start of Section \ref{partdif}, we showed that any biquotient of the form $M = N^m\times_{S^1} S^3$ has cohomology ring isomorphic to that of $N^m\times S^2$.  Thus, in order to apply Theorem \ref{finitediff}, we need only show all such biquotients have the same Pontrjagin classes.

\begin{proposition}\label{soclass}
For each $m$, the Pontrjagin classes of $M =N^m\times_{S^1} S^3$ are independent of the choice of $S^1$ action on $N^m$.

\end{proposition}

\begin{proof}

Recall that we have a bundle $N^m\rightarrow M\stackrel{\pi}{\rightarrow} S^2$.  Over a chart $\mathbb{R}^2\cong U\subseteq S^2$, the bundle trivializes, so $\pi^{-1}(U)\cong N^m\times \mathbb{R}^2$.  In particular, the normal bundle of $N^m$ in $M$ is trivial.

Now, let $i:N^m\rightarrow M$ be the inclusion of the fiber.  From the discussion at the beginning of this section, we know that $i^\ast$ is an isomorphism on cohomology groups in degree a multiple of $4$.  Since the normal bundle is trivial, it follows that $i^\ast$ identifies the Pontryagin classes of $N^m$ with those of $M$.

\end{proof}

This completes the proof of Theorem \ref{cclass}.


\end{document}